\documentclass[12pt]{article}
\usepackage{graphicx}
\usepackage{graphics}
\usepackage{amsthm}
\usepackage{amssymb, amsmath}
\title{Functional equations for the Stieltjes constants}  
\author{Mark W. Coffey\\
Department of Physics\\
Colorado School of Mines\\
Golden, CO  80401\\
mcoffey@mines.edu\\
(Received $\mbox{~~~~~~~~~~~~~~~~~~~~~~~~~~~~~~~2014}$)}
\date{February 23, 2014}
\begin{document}
\maketitle
\baselineskip=25 pt
\begin{abstract}

The Stieltjes constants $\gamma_k(a)$ appear as the coefficients in the regular
part of the Laurent expansion of the Hurwitz zeta function $\zeta(s,a)$ about $s=1$.
We present the evaluation of $\gamma_1(a)$ and $\gamma_2(a)$ at rational argument, being 
of interest to theoretical and computational analytic number theory and elsewhere.  
We give multiplication formulas for $\gamma_0(a)$, $\gamma_1(a)$, and $\gamma_2(a)$, and
point out that these formulas are cases of an addition formula previously presented.
We present certain integral evaluations generalizing Gauss' formula for the 
digamma function at rational argument.  In addition, we give the asymptotic form of
$\gamma_k(a)$ as $a \to 0$ as well as a novel technique for evaluating integrals with
integrands with $\ln(-\ln x)$ and rational factors.


\end{abstract}
 
\vspace{.25cm}
\baselineskip=15pt
\centerline{\bf Key words and phrases}
\medskip 

\noindent

Stieltjes constants, Riemann zeta function, Hurwitz zeta function, Laurent expansion,
digamma function, polygamma function 

\vfill
\centerline{\bf 2010 AMS codes}   
11M35, 11M06, 11Y60.  secondary:  05A10

\baselineskip=25pt
\pagebreak
\medskip
\centerline{\bf Introduction and statement of results}
\medskip

The Stieltjes (or generalized Euler) constants $\gamma_k(a)$ appear as 
expansion coefficients in the Laurent series for the Hurwitz zeta function 
$\zeta(s,a)$ about its simple pole at $s=1$ 
\cite{briggs,coffeyjmaa,coffeydec11,hardy,kluyver,mitrovic,stieltjes,wilton2},
$$\zeta(s,a)={1 \over {s-1}}+\sum_{n=0}^\infty {{(-1)^n} \over {n!}}\gamma_n(a)(s-1)^n.
\eqno(1.1)$$
These constants are important in analytic number theory and elsewhere, where
they appear in various estimations and as a result of asymptotic analyses, being given
by the limit relation
$$\gamma_k(a)=\lim_{N \to \infty} \left[\sum_{j=1}^N {{\ln^k (j+a)} \over j}-{{\ln^{k+1}(N+a)} \over {k+1}}\right].$$
In particular, $\gamma_0(a)=-\psi(a)$, where $\psi(z)=\Gamma'(z)/\Gamma(z)$ is the
digamma function, with $\Gamma(z)$ the Gamma function.
With $\gamma$ the Euler constant and $\gamma_1=\gamma_1(1)$ and $\gamma_2=\gamma_2(1)$, we
recall the connection with sums of reciprocal powers of the nontrivial zeros $\rho$ of the
Riemann zeta function,
$$\sum_\rho {1 \over \rho^2}=1-{\pi^2 \over 8}+2\gamma_1+\gamma^2, ~~~~
\sum_\rho {1 \over \rho^3}=1-{7 \over 8}\zeta(3)+\gamma^3+3\gamma \gamma_1+{3 \over 2}\gamma_2,$$
such relations following from the Hadamard factorization.

An effective asymptotic expression for $\gamma_k$ \cite{knessl1} and $\gamma_k(a)$
\cite{knessl2} for $k \gg 1$ has recently been given.  From this expression, 
previously known results on sign changes within the sequence of Stieltjes constants follow.

In this paper, we first evaluate the first and second Stieltjes constants at rational
argument.  These decompositions are effectively Fourier series, thus implying many
extensions and applications, and they supplement the relations presented in 
\cite{coffeystdiffs}.  We then present multiplication formulas for the zeroth, first, and 
second Stieltjes constants, and certain log-log integrals with integer parameters.  The latter
integral evaluations also provide explicit expressions for the differences $\gamma_1(j/m)
-\gamma_1$ and $\gamma_2(j/m)-\gamma_2$.  Besides elaborating on a multiplication formula for
the Stieltjes constants, the Discussion section provides examples of integrals evaluating in
terms of differences of the first and second of these constants.  In addition, presented
there is a novel method of determining log-log integrals with a certain polynomial 
denominator integrand. 

We recall the connection of differences of Stieltjes constants with logarithmic sums,
$$\gamma_\ell(a)-\gamma_\ell(b)=\sum_{n=0}^\infty \left[{{\ln^\ell(n+a)} \over 
{n+a}}-{{\ln^\ell(n+b)} \over {n+b}}\right].  \eqno(1.2)$$
Very recently an evaluation of $\gamma_1(j/m)-\gamma_1$ has also appeared \cite{blago}.
However, the method of proof is circuitous--integrals are used in addition to multiple
applications of functional equations.  

The Hurwitz zeta function, initially defined by $\zeta(s,a)=\sum_{n=0}^\infty
(n+a)^{-s}$ for $\mbox{Re} ~s>1$, has an analytic continuation to the whole
complex plane \cite{berndt,fine,titch,karatsuba}.   
In the case of $a=1$, $\zeta(s,a)$ reduces to the Riemann zeta function
$\zeta(s)$ \cite{edwards,ivic,riemann}.
In this instance, by convention, the Stieltjes constants
$\gamma_k(1)$ are simply denoted $\gamma_k$ \cite{briggs,hardy,kluyver,kreminski,mitrovic,yue}.  We recall that
$\gamma_k(a+1)=\gamma_k(a)-(\ln^k a)/a$, and more generally that for $n \geq 1$ an
integer 
$$\gamma_k(a+n)=\gamma_k(a)-\sum_{j=0}^{n-1} {{\ln^k(a+j)} \over {a+j}},$$
as follows from the functional equation $\zeta(s,a+n)=\zeta(s,a)-\sum_{j=0}^{n-1}
(a+j)^{-s}$.  
In fact, an interval of length $1/2$ is sufficient to characterize the $\gamma_k(a)$'s \cite{hansen}.

Unless specified otherwise below, letters $j$, $k$, $\ell$, $m$, $n$, and $r$ denote positive integers.  The Euler constant is given by $\gamma=-\psi(1)=\gamma_0(1)$.  The polygamma
functions are denoted $\psi^{(n)}(z)$ and we note that $\psi^{(n)}(z)=(-1)^{n+1}n!\zeta(n+1,z)$
\cite{nbs,grad}.

{\bf Proposition 1}.  For $m>1$ and $j<m$, (a)
$$\gamma_1\left({j \over m}\right)=\gamma_1+\gamma^2+\gamma \ln 2\pi m+\ln (2\pi)\ln m+
{{\ln^2 m} \over 2}+(\gamma+\ln 2\pi m)\psi\left({j \over m}\right)$$
$$+\pi \sum_{r=1}^{m-1}\sin {{2\pi jr} \over m}\ln \Gamma \left({r \over m}\right)+
\sum_{r=1}^{m-1} \cos{{2\pi jr} \over m} \zeta''\left(0,{r \over m}\right), \eqno(1.3)$$
and (b)
$${1 \over 2}\gamma_2\left({j \over m}\right)={\gamma_2 \over 2}-\gamma_1\ln m+{\gamma \over
8}\pi^2+{\pi^2 \over 8}\ln m+{\gamma \over 2}\ln^2 m+{{\ln^3 m} \over 6}$$
$$+A+\sum_{r=1}^{m-1}\cos{{2\pi jr} \over m}\left\{\left[{\pi^2 \over 6}+(\gamma+\ln 2\pi m)^2
\right]\zeta'\left(0,{r \over m}\right)-(\gamma+\ln 2\pi m)\zeta''\left(0,{r \over m}\right)
+{1 \over 3}\zeta'''\left(0,{r \over m}\right)\right\}, \eqno(1.4)$$
where 
$$A={\pi^3 \over {24}}\cot {{\pi j} \over m}+{\pi^2 \over 8}\psi\left({j \over m}\right)$$
$$-{\pi \over 2}\sum_{r=1}^{m-1}\sin {{2\pi jr} \over m}\left\{-\zeta\left(0,{r \over m}
\right)\left[{\pi^2 \over 6}+(\gamma+\ln 2\pi m)^2\right]
+2(\gamma+\ln 2\pi m)\zeta'\left(0,{r \over m}\right)-\zeta''\left(0,{r \over m}
\right)\right\},$$
and $'$ indicates differentiation with respect to the first argument.  

Both parts (a) and (b) may be written in many alternative forms.  For instance, for (b),
the well known relation $\zeta'(0,a)=\ln \Gamma(a)-\ln(2\pi)/2$ may be used, and as well,
the right member may be modified by introducing $\gamma_1(j/m)$.

Various summation results are known for the Stieltjes constants, including \cite{coffeyjmaa}
$$\sum_{r=1}^q \gamma_k\left({r \over q}\right)=q(-1)^k {{\ln^{k+1}q} \over {k+1}}
+q\sum_{j=0}^k (-1)^j{k \choose j}(\ln^j q) \gamma_{k-j}.$$
As we briefly indicate, the $k=1$ and $2$ cases follow from Proposition 1.
{\newline \bf Corollary 1}.
$$\sum_{r=1}^q \gamma_1\left({r \over q}\right)=-{q \over 2}\ln^2 q+q(-\gamma \ln q+\gamma_1),$$
and
$$\sum_{r=1}^q \gamma_2\left({r \over q}\right)={q \over 3}\ln^3 q+q\left(\gamma \ln^2q -2\gamma_1\ln q+ \gamma_2\right).$$

{\it Proof}.  The summation for $\gamma_1$ follows from (1.3) as the sum over the cosine
and sine terms are zero, and we have the readily verified relation $\sum_{r=1}^q\psi\left(
{r \over q}\right)=-q(\gamma+\ln q)$.  Similarly for the summation for $\gamma_2$, the
sine and cosine terms do not contribute, the just-mentioned summation of $\psi(r/q)$ holds,   
and $\sum_{j=1}^{q-1} \cot {{\pi j} \over q}=0$. \qed

{\bf Corollary 2}.  For $\ell=1,2,\ldots,m-1$, (a)
$${1 \over 2}\sum_{j=1}^{m-1}\gamma_2\left({j \over m}\right)\sin {{2\pi j \ell} \over m}
={\pi^3 \over {48}}(2\ell-m)$$
$$-{{\pi m} \over 4}\left\{\left({\pi^2 \over 6}+(\gamma+\ln 2\pi m)^2
\right)\left[\zeta\left(0,1-{\ell \over m}\right)-\zeta\left(0,{\ell \over m}\right)\right]
\right.$$
$$\left.+2(\gamma+\ln 2\pi m)\left[\zeta'\left(0,{\ell \over m}\right)-\zeta'\left(0,1-{\ell \over m}\right)\right]+\zeta''\left(0,1-{\ell \over m}\right)-\zeta''\left(0,{\ell \over m}\right)\right\},$$
and (b)
$${1 \over 2}\sum_{j=1}^{m-1}\gamma_2\left({j \over m}\right)\cos {{2\pi j \ell} \over m}
=-\left[{\gamma_2 \over 2}-\gamma_1\ln m+{\gamma \over 8}\pi^2+{\pi^2 \over 8}\ln m+{\gamma
\over 2}\ln^2 m+{1 \over 6}\ln^3 m\right]$$
$$+{\pi^2 \over 8}\left[m\ln\left(2\sin {{\pi \ell} \over m}\right)+\gamma\right]$$
$$-\sum_{r=1}^{m-1}\left\{\left({\pi^2 \over 6}+(\gamma+\ln 2\pi m)^2\right)\zeta'\left(0,{r \over m}\right)-(\gamma+\ln 2\pi m)\zeta''\left(0,{r \over m}\right)+{1 \over 3}\zeta''' \left(0,{r \over m}\right)\right\}$$
$$+{m \over 2}\left\{\left({\pi^2 \over 6}+(\gamma+\ln 2\pi m)^2\right)\left[\zeta'\left(
0,{\ell \over m}\right)+\zeta'\left(0,1-{\ell \over m}\right)\right]\right.$$
$$\left.-(\gamma+\ln 2\pi m)\left[\zeta''\left(0,{\ell \over m}\right)+\zeta''\left(0,1-{\ell \over m}\right)\right]
+{1 \over 3}\left[\zeta'''\left(0,{\ell \over m}\right)+\zeta'''\left(0,1-{\ell \over m}\right)\right]\right\}.$$
Part (b) may be rewritten by using the three sums
$$\sum_{\ell=1}^{m-1}\zeta'\left(0,{\ell \over m}\right)=-{{\ln m} \over 2}, ~~~~~~
\sum_{\ell=1}^{m-1}\zeta''\left(0,{\ell \over m}\right)=-{1 \over 2}\ln^2 m-(\ln m)\ln (2\pi),$$
and 
$$\sum_{\ell=1}^{m-1}\zeta'''\left(0,{\ell \over m}\right)=-{{\ln^3 m} \over 2}-{3 \over 2}
\ln^2 m \ln(2\pi) +3(\ln m)\left({\gamma^2 \over 2}-{\pi^2 \over {24}}-{1 \over 2}\ln^2(2\pi)
+\gamma_1\right).$$

{\bf Corollary 3}.  (Asymptotic expressions)  For $m\to \infty$,
$$\gamma_1\left({1 \over m}\right) \sim -m\ln m+\gamma_1,$$
and
$$\gamma_2\left({1 \over m}\right) \sim m\ln^2m +\gamma_2.$$
The general situation for $\gamma_k(a)$ with $a\to 0$ is more conveniently proved
otherwise, and is presented in Proposition 5.

{\bf Proposition 2}.  
For 
Re $z>0$ and $0<k<2$,
$$\psi(kz)=\sum_{n=0}^\infty (k-1)^n {z^n \over {n!}}\psi^{(n)}(z).$$

{\bf Proposition 3}.
For 
Re $z>0$ and $0<k<2$, (a)
$$\gamma_1(kz)=\gamma_1(z)+\gamma [\psi(kz)-\psi(z)]-\sum_{n=1}^\infty (1-k)^n z^n\left[ {{(-1)^{n+1}} \over {n!}} \psi^{(n)}(z)\psi(n+1)+\left.{{\partial \zeta(s,z)} \over {\partial s}}\right|_{s=n+1}\right],$$
and (b)
$$\gamma_2(kz)=\gamma_2(z)+[\gamma^2-\zeta(2)][\psi(z)-\psi(kz)]$$
$$+2\sum_{n=1}^\infty (1-k)^n z^n\left[\left(\gamma+{1 \over 2}\psi(n+1)\right)\psi(n+1)
\zeta(n+1,z)+{1 \over 2}\psi'(n+1)\zeta(n+1,z) \right.$$
$$\left. +[\gamma+\psi(n+1)]\left.{{\partial \zeta(s,z)} \over {\partial s}}\right|_{s=n+1}+{1 \over 2}\left.{{\partial^2 \zeta(s,z)} \over {\partial s^2}}\right|_{s=n+1}\right].$$

The following result, wherein the differences $\gamma_1(j/m)-\gamma_1$ and $\gamma_2(j/m)-\gamma_2$ appear, is a generalization of Gauss' formula for the digamma function at rational argument.  For integers $p$ and $q$ with $0<p<q$, we put
$$I^k_{pq}\equiv q\int_0^1 \left({{x^{q-1}-x^{p-1}} \over {1-x^q}}\right)\ln^k(-\ln x)dx.$$
{\bf Proposition 4}.  Let $\omega_k\equiv \exp(2\pi ik/q)$.  Then (a)
$$I_{pq}^1=-\sum_{k=1}^{q-1}(\omega_k^p-1)\left[\gamma\ln\left({{\omega_k-1} \over \omega_k}
\right)+\sum_{n=1}^\infty \zeta'(1-n){{\ln^n(\omega_k^{-1})} \over {n!}}-\gamma \ln(-\ln
\omega_k^{-1})-{1 \over 2}\ln^2(-\ln \omega_k^{-1})\right]$$
$$-q\left(\gamma_1+{\gamma^2 \over 2}+{\pi^2 \over {12}}\right)$$
$$=-(\gamma+\ln q)\left[\gamma+\psi\left({p \over q}\right)\right]+\gamma_1\left({p \over q}
\right)-\gamma_1,$$
and (b)
$$I_{pq}^2=\sum_{k=1}^{q-1}(\omega_k^p-1)\left[\left(\gamma^2+{\pi^2 \over 6}\right)
\ln\left({{\omega_k-1} \over \omega_k}\right)+2\gamma \left.{{\partial \mbox{Li}_s} \over {\partial s}}(\omega_k^{-1})\right|_{s=1}-\left.{{\partial^2 \mbox{Li}_s} \over {\partial s^2}}(\omega_k^{-1})\right|_{s=1}\right]$$
$$=\left(\gamma^2+{\pi^2 \over 6}+2\gamma \ln q+\ln^2 q\right)\left[\gamma+\psi\left({p \over q}\right)\right]
-2(\gamma+\ln q)\left[\gamma_1\left({p \over q}\right)-\gamma_1\right]+\gamma_2- \gamma_2\left({p \over q}\right).$$
In (b), the analytically continuable polylogarithm (or Jonqui\`{e}re) function
Li$_s(z)=\sum_{k=1}^\infty z^k/k^s$.  Explicit expressions for the partial derivatives
appearing there are provided.

{\bf Proposition 5}.  (Asymptotic expression) For $a \to 0$,
$$\gamma_k(a) \sim {1 \over a}\ln^k a + \gamma_k.$$

\bigskip
\centerline{\bf Proof of Propositions}
\medskip

{\it Proposition 1}. We will be expanding a functional equation due to Hurwitz
(\cite{apostol}, p. 261; \cite{hurwitz}, p. 93),
$$\zeta\left(s,{j \over m}\right)={{2^s \pi^{s-1}} \over m^{1-s}}\Gamma(1-s)\sum_{r=1}^m
\sin\left({{\pi s} \over 2}+{{2\pi j r} \over m}\right)\zeta\left(1-s,{r \over m}\right),
\eqno(2.1)$$
holding for $1 \leq j \leq m$, about $s=1$.  We will use two elementary trigonometric
identities
$$-2\sum_{r=1}^m \cos {{2\pi j r} \over m}\zeta\left(0,{r \over m}\right)
=-2\sum_{r=1}^m \left({1 \over 2}-{r \over m}\right)\cos {{2\pi j r} \over m}=1,
~~~~1 \leq j<m, \eqno(2.2a)$$
and
$$\sum_{r=1}^m \sin {{2\pi j r} \over m}\zeta\left(0,{r \over m}\right)
=\sum_{r=1}^m \left({1 \over 2}-{r \over m}\right)\sin {{2\pi j r} \over m}={1 \over 2}
\cot{{\pi j} \over m}, ~~~~1 \leq j<m, \eqno(2.2b)$$
as well as the values
$$\zeta''(0)=\gamma_1+{\gamma^2 \over 2}-{\pi^2 \over {24}}-{{\ln^2 2\pi} \over 2}, \eqno(2.3a)$$
and
$$\zeta'''(0)=\gamma^3+{3 \over 2}\gamma^2 \ln 2\pi-{\pi^2 \over 8}\ln 2\pi-{1 \over 2}
\ln^3 2\pi+3(\gamma+\ln 2\pi) \gamma_1+{3 \over 2}\gamma_2-\zeta(3). \eqno(2.3b)$$
These derivative values may be obtained via (1.1) and the functional equation of the zeta
function, $\zeta(s)=2(2\pi)^{s-1} \sin(\pi s/2)\Gamma(1-s)\zeta(1-s)$.
For (2.2) we have used the well known relation $\zeta(0,a)=1/2-a=-B_1(a)$, where $B_1(a)$
is the first-degree Bernoulli polynomial.  We recall that about $s=1$, $\Gamma(1-s)$ has a
simple pole and that
$$\Gamma(1-s)=-{1 \over {s-1}}-\gamma-\left({\pi^2 \over {12}}+{\gamma^2 \over 2}\right)
(s-1)-{1 \over 6}\left(\gamma^3+{\gamma \over 2}\pi^2-\psi''(1)\right)(s-1)^2+O[(s-1)^3],$$
wherein the tetragamma function value $\psi''(1)=-2\zeta(3)$, and recall the expansion
$(2\pi m)^{s-1}=\sum_{\ell=0}^\infty {{\ln^\ell 2\pi m} \over {\ell!}}(s-1)^\ell$.
For the sine factor on the right side of (2.1) we have
$$\sin\left({{\pi s} \over 2}+{{2\pi j r} \over m}\right)=\cos{{2\pi jr} \over m}
\sum_{\ell=0}^\infty {{(-1)^\ell \pi^{2\ell}} \over {(2\ell)!2^{2\ell}}}(s-1)^{2\ell}
+\sin{{2\pi jr} \over m}\sum_{\ell=0}^\infty {{(-1)^{\ell+1} \pi^{2\ell+1}} \over {(2\ell+1)!2^{2\ell+1}}}(s-1)^{2\ell+1}.$$

The left side of (2.1) expands as
$$\zeta\left(s,{j \over m}\right)={1 \over {s-1}}-\psi\left({j \over m}\right)-\gamma_1
\left({j \over m}\right)(s-1)+{1 \over 2}\gamma_2\left({j \over m}\right)(s-1)^2 +
O[(s-1)^3].$$
The polar contributions on the two sides of (2.1) cancel due to the cosine sum (2.2a).
At order $(s-1)^0$, one finds
$$\psi\left({j \over m}\right)=\sum_{r=1}^m\left\{-\pi \zeta\left(0,{r \over m}\right)
\sin\left({{2\pi jr} \over m}\right)+2\cos{{2\pi jr} \over m}\left[\zeta\left(0,{r \over m}\right)(\gamma+\ln 2\pi m)-\zeta'\left(0,{r \over m}\right)\right]\right\}, \eqno(2.4)$$
where both of the sums of (2.2) apply.  With $\zeta'(0,r/m)=\ln \Gamma(r/m)-\ln(2\pi)/2$,
$$\psi\left({j \over m}\right)=-\gamma-\ln 2\pi m-{\pi \over 2}\cot {{\pi j} \over m}
-2\sum_{r=1}^{m-1} \cos{{2\pi jr} \over m}\ln \Gamma\left({r \over m}\right),$$
being a form of one of Gauss' formulas for $\psi$ at rational argument.

At order $(s-1)^1$, one finds
$$-\gamma_1\left({j \over m}\right)=\sum_{r=1}^m\left\{ {\pi^2 \over 4}\cos{{2\pi jr} \over m}\zeta\left(0,{r \over m}\right)
-{\pi \over 2}\sin{{2\pi jr} \over m}\left[\zeta\left(0,{r \over m}\right)(\gamma+\ln 2\pi m)
+2\zeta'\left(0,{r \over m}\right)\right] \right.$$
$$\left. +\cos{{2\pi jr} \over m}\left[-\zeta\left(0,{r \over m}\right)\left({\pi^2 \over 6}
+(\gamma+\ln 2\pi m)^2\right)+2(\gamma+\ln 2\pi m)\zeta'\left(0,{r \over m}\right)-
\zeta''\left(0,{r \over m}\right)\right]\right\}, \eqno(2.5)$$
and again the trigonometric summations of (2.2) apply.  
We add and subtract $(\gamma+\ln 2\pi m)\cos(2\pi jr/m)\zeta(0,r/m)$ to introduce the $\psi(j/m)$ expression (2.4).  Then separating the $r=m$ terms of the sums and using $\zeta''(0)$ from (2.3a) gives part (a).  

For part (b), at order $(s-1)^2$ in (2.1) we have
$${1 \over 2}\gamma_2\left({j \over m}\right)=A+B+C,$$
where
$$A=\sum_{r=1}^m\left\{-{\pi^3 \over {24}}\zeta\left(0,{r \over m}\right)\sin {{2\pi jr} \over
m}+{\pi^2 \over 4}\cos {{2\pi jr} \over m}\left[\zeta\left(0,{r \over m}\right)(\gamma+
\ln 2\pi m)-\zeta'\left(0,{r \over m}\right)\right]\right.$$
$$\left.-{\pi \over 2}\sin{{2\pi jr} \over m}\left[-\zeta\left(0,{r \over m}\right)\left({\pi^2 \over 6}+(\gamma+\ln 2\pi m)^2\right)+2(\gamma+\ln 2\pi m)\zeta'\left(0,{r \over m}\right)
-\zeta''\left(0,{r \over m}\right)\right ]\right \},$$
$$B=-\sum_{r=1}^m \cos {{2\pi jr} \over m}\zeta\left(0,{r \over m}\right)\left[{1 \over 3}
\ln^3 2\pi m-\left(\gamma^2+{\pi^2 \over 6}\right)\ln 2\pi m-\gamma\ln^2 2\pi m \right.$$
$$\left. +{1 \over 3} \left(\gamma^3-{\gamma \over 2}\pi^2+2\zeta(3)\right)\right],$$
and
$$C=\sum_{r=1}^m \cos {{2\pi jr} \over m}\left\{\left[{\pi^2 \over 6}+(\gamma+\ln 2\pi m)^2
\right]\zeta'\left(0,{r \over m}\right)-(\gamma+\ln 2\pi m)\zeta''\left(0,{r \over m}\right)
\right.$$
$$\left. +{1 \over 3}\zeta'''\left(0,{r \over m}\right) \right\}.$$
The sum in $B$ is immediately evaluated with the aid of (2.2a).  The expression for $A$ is
rewritten by using (2.4), making use of the values $\psi(j/m)$, and (2.2b).  In the sum $C$, 
the $r=m$ terms are separated and all of the values $\zeta'(0)=-\ln(2\pi)/2$, (2.3a), and (2.3b) are used.  Combining terms then yields (1.4).  \qed

{\it Remarks}.  The expression for $C$ in the proof may be written in terms of $\gamma_1(j/m)$
by using (2.5). 
For the evaluation of general $\gamma_k(j/m)$, there will always be a sum
$\sum_{r=1}^{m-1} \cos{{2\pi jr} \over m}\zeta^{(k+1)}\left(0,{r \over m}\right)$,
as according to the highest derivative of $\zeta(s,a)$ present, there will be no
derivative of the sine factor in (2.1).  The evaluation will then also contain a 
contribution of $\zeta^{(k+1)}(0,1)=\zeta^{(k+1)}(0)$.

{\it Corollary 2}.  
The proof uses the discrete orthogonality of sine and cosine functions, implying
$$\sum_{j=1}^{m-1}\cot{{\pi j} \over m}\cos{{2\pi j\ell} \over m}=0,$$
and
$$\sum_{j=1}^{m-1}\cot{{\pi j} \over m}\sin{{2\pi j\ell} \over m}=m-2\ell,$$
along with
$$\sum_{j=1}^{m-1}\psi\left({j \over m}\right)\cos {{2\pi j \ell} \over m}=m\ln\left(2\sin
{{\pi \ell} \over m}\right)+\gamma,$$
and
$$\sum_{j=1}^{m-1}\psi\left({j \over m}\right)\sin {{2\pi j \ell} \over m}={\pi \over 2}
(2\ell-m).$$
The three sums following part (b) may be determined by successively differentiating the
relation $\sum_{\ell=1}^{m-1}\zeta\left(s,{\ell \over m}\right)=(m^s-1)\zeta(s)$ and putting
$s$ to $0$. \qed

{\it Proposition 2}.  We will indicate four proofs.  For the first, we may use a standard
integral representation for the polygamma function \cite{grad} (p. 943) to write
$$\sum_{n=0}^\infty (k-1)^n {z^n \over {n!}}\psi^{(n)}(z)=\psi(z)+\sum_{n=1}^\infty (k-1)^n {z^n
\over {n!}}\int_0^1 {{t^{z-1}\ln^n t} \over {t-1}}dt$$
$$=\psi(z)+\int_0^1 {{[1-t^{(k-1)z}]} \over {t-1}}t^{z-1}dt$$
$$=\psi(z)-\psi(z)+\psi(kz)=\psi(kz).$$
\qed

For the second, we may start with the expansion \cite{grad} (p. 944)
$$\psi(x+1)=\psi(x)+{1 \over x}=-\gamma+\sum_{j=2}^\infty (-1)^j \zeta(j)x^{j-1}, ~~~~|x|<1.$$
Then
$$\psi^{(n)}(z)={{(-1)^{n+1}n!} \over z^{n+1}}+\sum_{j=2}^\infty (-1)^j \zeta(j){{(j-1)!} \over
{(j-n-1)!}}z^{j-n-1},$$
so that
$$\sum_{n=0}^\infty (k-1)^n {z^n \over {n!}}\psi^{(n)}(z)=\psi(z)+\sum_{n=1}^\infty (k-1)^n 
{z^n \over {n!}}\left[{{(-1)^{n+1}n!} \over z^{n+1}}+\sum_{j=2}^\infty (-1)^j \zeta(j){{(j-1)!} \over {(j-n-1)!}}z^{j-n-1}\right]$$
$$=\psi(z)-{1 \over {kz}}+{1 \over z}+\sum_{j=2}^\infty (-1)^j \zeta(j)z^{j-1}\sum_{n=1}^{j-1}(k-1)^n {{j-1} \choose n}$$
$$=\psi(z)-{1 \over {kz}}+{1 \over z}+\sum_{j=2}^\infty (-1)^j \zeta(j)z^{j-1}(k^{j-1}-1)$$
$$=\psi(z)-{1 \over {kz}}+{1 \over z}+\psi(kz)+{1 \over {kz}}-\psi(z)-{1 \over z}=\psi(kz).$$
\qed

A third method of proof follows from the representations
$$\gamma_0(a)=-\psi(a)=-\ln a-\sum_{k=1}^\infty {1 \over {k+1}}\sum_{\ell=0}^k (-1)^\ell{k \choose \ell} \ln(\ell+a)$$
and
$$\psi^{(j)}(a)=(-1)^{j-1}{{(j-1)!} \over a^j}+(j-1)!\sum_{k=1}^\infty {1 \over {k+1}}\sum_{\ell=0}^k (-1)^\ell{k \choose \ell} {{(-1)^{j-1}} \over {(\ell+a)^j}}. $$
\qed

The fourth method of proof follows that of the first proof of Proposition 3 so we omit it.

{\it Proposition 3}.  We have the multiplication formula
$$\zeta(s,kz)=\sum_{n=0}^\infty {{s+n-1} \choose n}(1-k)^n z^n \zeta(s+n,z),$$
being a case of a more general result of Truesdell \cite{truesdell}.
We then expand both sides about $s=1$.  Equating the coefficients of $(s-1)^0$ on both
sides gives another means of demonstrating Proposition 2.  Equating the coefficients of
$(s-1)^1$ and using Proposition 2 gives part (a).  Equating the coefficients of $(s-1)^2$
gives part (b).  \qed

For a second method of proof of part (a), we may use the integral representation for Re $s>1$ and Re $a>0$,
$$\zeta(s,a)={1 \over {\Gamma(s)}}\int_0^\infty {{t^{s-1}e^{-(a-1)t}} \over {e^t-1}}dt,$$
so that
$${{\partial \zeta(s,a)} \over {\partial s}}=-\psi(s)\zeta(s,a)+{1 \over {\Gamma(s)}}\int_0^\infty {{t^{s-1}e^{-(a-1)t}} \over {e^t-1}}\ln t ~dt.$$
Then
$$\sum_{n=1}^\infty (1-k)^n z^n \left.{{\partial \zeta(s,z)} \over {\partial s}}\right|_{s=n+1}
=\sum_{n=1}^\infty (1-k)^n z^n \left[-\psi(n+1)\zeta(n+1,z)+{1 \over {n!}}\int_0^\infty {{t^n e^{-(a-1)t}} \over {e^t-1}}\ln t ~dt\right]$$
$$=-\sum_{n=1}^\infty (1-k)^n z^n \psi(n+1)\zeta(n+1,z)+\int_0^\infty e^{-(z-1)t}{{[e^{-(k-1)tz}
-1]} \over {e^t-1}}\ln t ~dt.$$
Now
$$\int_0^\infty t^{s-1}e^{-(z-1)t}{{[e^{-(k-1)tz}-1]} \over {e^t-1}}dt=\Gamma(s)[\zeta(s,kz)-\zeta(s,z)],$$
giving
$$\int_0^\infty t^{s-1}e^{-(z-1)t}{{[e^{-(k-1)tz}-1]} \over {e^t-1}}\ln t ~dt=\Gamma(s)\psi(s)[\zeta(s,kz)-\zeta(s,z)]+\Gamma(s)\left[{{\partial \zeta}
\over {\partial s}}(s,kz)-{{\partial \zeta}\over {\partial s}}(s,z)\right].$$
Since by (1.1)
$${{\partial \zeta}\over {\partial s}}(s,a)=-{1 \over {(s-1)^2}}+\sum_{n=0}^\infty {{(-1)^{n+1}}
\over {n!}}\gamma_{n+1}(a)(s-1)^n$$
$$=-{1 \over {(s-1)^2}}-\gamma_1(a)(s-1)^0+{{\gamma_2(a)} \over 2}(s-1)^1 +O[(s-1)^2],$$
we find
$$\int_0^\infty e^{-(z-1)t}{{[e^{-(k-1)tz}-1]} \over {e^t-1}}\ln t ~dt=\gamma[\psi(kz)-\psi(z)]
+\gamma_1(kz)-\gamma_1(z),$$
and the Proposition again follows.  \qed

{\it Remarks}.  Obviously we may evaluate the integrals
$$\int_0^\infty e^{-(z-1)t}{{[e^{-(k-1)tz}-1]} \over {e^t-1}}\ln^j t ~dt$$
in terms of the difference $\gamma_j(kz)-\gamma_j(z)$.

The harmonic numbers $H_n=\sum_{k=1}^n 1/k=\psi(n)+\gamma$ and generalized harmonic
numbers 
$$H_n^{(r)}=\sum_{k=1}^n {1 \over k^r} ={{(-1)^{r-1}} \over {(r-1)!}}\left[\psi^{(r-1)}(n+1)
-\psi^{(r-1)}(1)\right]={{(-1)^{r-1}} \over {(r-1)!}}\int_0^1{{(t^n-1)} \over {t-1}}\ln^{r-1}t
~dt$$  
enter the representations of Proposition 3 and for the higher Stieltjes constants.
This is part of the elaboration of the following discussion section.

{\it Proposition 4}.  We first write, for $0<p<q$,
$$I^k_{pq}\equiv q\int_0^1 \left({{x^{q-1}-x^{p-1}} \over {1-x^q}}\right)\ln^k(-\ln x)dx
=\sum_{k=1}^{q-1} \int_0^1 \left({{\omega_k^p-1} \over {x-\omega_k}}\right)\ln^k(-\ln x)dx.$$
By performing logarithmic differentiation on the integral
$$\int_0^1 {{(-\ln x)^a} \over {x-\omega}}dx=-\Gamma(a)\mbox{Li}_a(\omega^{-1}),$$
one then has the following expressions:
$$I_{pq}^1=-\sum_{k=1}^{q-1}(\omega_k^p-1)\left[\gamma\ln\left({{\omega_k-1} \over \omega_k}
\right)+\left.{{\partial \mbox{Li}_s} \over {\partial s}}(\omega_k^{-1})\right|_{s=1}\right]$$
and
$$I_{pq}^2=\sum_{k=1}^{q-1}(\omega_k^p-1)\left[\left(\gamma^2+{\pi^2 \over 6}\right)
\ln\left({{\omega_k-1} \over \omega_k}\right)+2\gamma \left.{{\partial \mbox{Li}_s} \over {\partial s}}(\omega_k^{-1})\right|_{s=1}-\left.{{\partial^2 \mbox{Li}_s} \over {\partial s^2}}(\omega_k^{-1})\right|_{s=1}\right].$$

We next present the partial derivatives of the polylogarithm function.  These result from
expansion of the following expression in powers of $s-1$ \cite{erdelyi}:
$$\mbox{Li}_s(z)=\Gamma(1-s)\ln^{s-1}\left({1 \over z}\right)+\sum_{n=0}^\infty \zeta(s-n)
{{\ln^n z} \over {n!}}, ~~~~|\ln z|<2\pi,$$
wherein the polar part of the first term on the right is cancelled by the pole $1/(s-1)$
of the $n=0$ term of the sum.  We obtain:
$$\left.{{\partial \mbox{Li}_s} \over {\partial s}}(z)\right|_{s=1}=-\gamma_1-{\gamma^2 \over
2}-{\pi^2 \over {12}}-\gamma\ln(-\ln z)-{1 \over 2}\ln^2(-\ln z)
+\sum_{n=1}^\infty \zeta'(1-n){{\ln^n z} \over {n!}}, \eqno(2.6)$$
and
$$\left.{{\partial^2 \mbox{Li}_s} \over {\partial s^2}}(z)\right|_{s=1}={1 \over 6}\left[
-2\gamma^3-\gamma \pi^2-6\gamma^2\ln(-\ln z)-\pi^2\ln(-\ln z)-6\gamma \ln^2(-\ln z) \right.$$
$$\left. -2\ln^3(-\ln z) -4\zeta(3)\right]
+\gamma_2+\sum_{n=1}^\infty \zeta''(1-n){{\ln^n z} \over {n!}}. \eqno(2.7)$$
Part (a) then makes use of the first derivative and the sum $\sum_{k=1}^{q-1}(\omega_k^p-1)=-q$.

For the second evaluation of $I_{pq}^k$ we use
$$I_{pq}^k=q\sum_{m=0}^\infty \int_0^1 (x^{q-1}-x^{p-1})x^{qm}\ln^k(-\ln x)dx$$
$$=q\sum_{m=0}^\infty \int_0^\infty \left(e^{-(q-1)u}-e^{-(p-1)u}\right)e^{-(qm+1)u}\ln^k u ~du.$$
By using logarithmic differentiation of the Gamma function integral,
$$I_{pq}^1=q\sum_{m=0}^\infty \left[{{\gamma+\ln(p+mq)} \over {p+mq}}-{{\gamma+\ln(m+1)q} \over {(m+1)q}}\right]$$
$$=-(\gamma+\ln q)\left[\gamma+\psi\left({p \over q}\right)\right]+\sum_{m=0}^\infty \left[
{{\ln(m+p/q)} \over {m+p/q}}-{{\ln(m+1)} \over {m+1}}\right]$$
$$=-(\gamma+\ln q)\left[\gamma+\psi\left({p \over q}\right)\right]+\gamma_1\left({p \over q}
\right)-\gamma_1,$$
wherein we applied (1.2) and the well known summation (e.g., \cite{grad}, p. 943)
$$\psi(x)=-\gamma-\sum_{k=0}^\infty \left({1 \over {x+k}}-{1 \over {k+1}}\right).$$

The other evaluation of (b) goes similarly, with
$$I_{pq}^2=q\sum_{m=0}^\infty \int_0^\infty \left(e^{-(q-1)u}-e^{-(p-1)u}\right)e^{-(qm+1)u}
\ln^2 u ~du$$
$$=\sum_{m=0}^\infty \left\{{1 \over {6(m+1)}}\left[6\gamma^2+\pi^2+12\gamma \ln(m+1)q+6\ln^2
(m+1)q\right]\right.$$
$$\left.-{1 \over {6(m+p/q)}}\left[6\gamma^2+\pi^2+12\gamma \ln(mq+p)q+6\ln^2
(mq+p)\right]\right\}$$
$$=\left(\gamma^2+{\pi^2 \over 6}\right)\left[\gamma+\psi\left({p \over q}\right)\right]
+2\gamma\sum_{m=0}^\infty\left[{{\ln(m+1)q} \over {m+1}}-{{\ln(mq+p)} \over {m+p/q}}\right]$$
$$+\sum_{m=0}^\infty\left[{{\ln^2(m+1)q} \over {m+1}}-{{\ln^2(mq+p)} \over {m+p/q}}\right]$$
$$=\left(\gamma^2+{\pi^2 \over 6}+2\gamma \ln q\right)\left[\gamma+\psi\left({p \over q}\right)\right]-2\gamma\left[\gamma_1\left({p \over q}\right)-\gamma_1\right]$$
$$+\sum_{m=0}^\infty\left[{{\ln^2(m+1)} \over {m+1}}-{{\ln^2(m+p/q)} \over {m+p/q}}
+2\ln q\left({{\ln(m+1)} \over {m+1}}-{{\ln(m+p/q)} \over {m+p/q}}\right)\right.$$
$$\left. +\ln^2q \left({1 \over {m+1}}-{1 \over {m+p/q}}\right)\right]$$
$$=\left(\gamma^2+{\pi^2 \over 6}+2\gamma \ln q+\ln^2 q\right)\left[\gamma+\psi\left({p \over q}\right)\right]
-2(\gamma+\ln q)\left[\gamma_1\left({p \over q}\right)-\gamma_1\right]+\gamma_2- \gamma_2\left({p \over q}\right).$$
\qed

{\it Remark}.  For applications or computation with Proposition 4, it is important that
the values of $\ln(\pm \omega_k^{\pm 1})$ are kept to the principal branch, e.g., with 
$-\pi <\mbox{Im} ~\ln z \leq \pi$.  This requirement maintains a real-valued result for
$I_{pq}^k$.

{\it Elaboration of the partial derivative (2.6)}.

\baselineskip=15pt
The partial derivative (2.6) may also be represented as
$$-\left.{{\partial \mbox{Li}_s(z)} \over {\partial s}}\right|_{s=1}=\int_1^\infty {z^x \over x}\ln x ~dx +\int_1^\infty z^x\left({{(\ln z)\ln x} \over x}+{1 \over x^2}-{{\ln x} \over x^2} 
\right)P_1(x)dx, \eqno(2.8)$$
following from \cite{coffeyaddison}
$$-{{\partial \mbox{Li}_s(z)} \over {\partial s}}=\int_1^\infty {z^x \over x^s}
\ln x ~dx +\int_1^\infty z^x\left({{(\ln z)\ln x} \over x^s}+{1 \over x^{s+1}}-s {{\ln x} \over x^{s+1}}\right)P_1(x)dx,$$
wherein $P_1(x) \equiv B_1(x-[x])$.
We relate (2.6) and (2.8).  By using the expansion $z^x=\sum_{j=0}^\infty \ln^j z{x^j \over {j!}}$, we have
$$\int_1^a {z^x \over x}\ln x ~dx=\sum_{j=0}^\infty {{\ln^j z} \over {j! j^2}}[1+a^j(j\ln a-1)]$$
$$=-\gamma \ln a-\Gamma(0,-a\ln z)\ln a+{{\ln^2 a} \over 2}+\ln z ~_3F_3(1,1,1;2,2,2;\ln z)$$
$$-a\ln z ~_3F_3(1,1,1;2,2,2;a\ln z)-\ln a \ln(-a\ln z),$$
where $_pF_q$ is the generalized hypergeometric function and the incomplete Gamma function
$\Gamma(\alpha,x)=\Gamma(\alpha)-\sum_{n=0}^\infty {{(-1)^nx^{n+\alpha}} \over {n!(n+\alpha)}}$.  
The asymptotic form of the $_3F_3$ function as $a \to \infty$ may be considered as in \cite{coffeyhypergeo} and the result is
$$\int_1^\infty  {z^x \over x}\ln x ~dx={\gamma^2 \over 2}+{\pi^2 \over {12}}+\ln z 
~_3F_3(1,1,1;2,2,2;\ln z)+\gamma \ln (-\ln z)+{1 \over 2}\ln^2(-\ln z), ~~~~|z|<1. \eqno(2.9)$$
This result (2.9) may also be obtained as a reduction of a Meijer-$G$ function.  However,
we omit details of this evaluation.

By using respectively the partial derivative of (1.1) with respect to $s$ and then (2.8) and (2.9) we have these additional expressions for the partial derivative (2.6):
$$\left.{{\partial \mbox{Li}_s} \over {\partial s}}(z)\right|_{s=1}=-\gamma_1-{\gamma^2 \over
2}-{\pi^2 \over {12}}-\gamma\ln(-\ln z)-{1 \over 2}\ln^2(-\ln z)$$
$$-\ln z ~_3F_3(1,1,1;2,2,2;\ln z)-\sum_{j=0}^\infty {\gamma_{j+1} \over {j!}}\sum_{n=1}^\infty
{{\ln^n z} \over {n!}}n^j$$
$$=-{\gamma^2 \over 2}-{\pi^2 \over {12}}-\gamma\ln(-\ln z)-{1 \over 2}\ln^2(-\ln z)
-\ln z ~_3F_3(1,1,1;2,2,2;\ln z)$$
$$-\int_1^\infty z^x\left({{(\ln z)\ln x} \over x}+{1 \over x^2}-{{\ln x} \over x^2}\right)
P_1(x)dx$$
$$=-{\gamma^2 \over 2}-{\pi^2 \over {12}}-\gamma\ln(-\ln z)-{1 \over 2}\ln^2(-\ln z)
-\ln z ~_3F_3(1,1,1;2,2,2;\ln z)$$
$$-\sum_{j=1}^\infty \int_0^1 z^{y+j}\left[{{(\ln z)\ln(y+j)} \over {y+j}}+{1 \over {(y+j)^2}}
-{{\ln(y+j)} \over {(y+j)^2}}\right]\left(y-{1 \over 2}\right)dy.$$
By comparing (2.6) with the second expression above we conclude that
$$\gamma_1-\sum_{n=1}^\infty \zeta'(1-n){{\ln^n z} \over {n!}}$$
$$=\ln z ~_3F_3(1,1,1;2,2,2;\ln z)
+\int_1^\infty z^x\left({{(\ln z)\ln x} \over x}+{1 \over x^2}-{{\ln x} \over x^2}\right)
P_1(x)dx$$

\baselineskip=25pt
{\it Proposition 5}.  Let $C_k(a) \equiv \gamma_k(a) -(\ln^k a)/a$.  With $B_n(x)$ the
Bernoulli polynomials, their periodic extension is denoted $P_n(x) \equiv B_n(x-[x])$
and we have the representation \cite{yue}
$$C_n(a) = (-1)^{n-1} n!\sum_{k=0}^{n+1} {{s(n+1,n+1-k)} \over {k!}}
\int_1^\infty P_n(x-a) {{\ln^k x} \over x^{n+1}} dx, ~~n \geq 1, \eqno(2.10)$$
with $s(n,k)$ the Stirling numbers of the first kind.
We recall the Fourier expansions of $P_n(x)$ \cite{nbs} (p. 805),
$$P_{2n}(x-a)=(-1)^{n-1} {{2(2n)!} \over {(2\pi)^{2n}}}\sum_{k=1}^\infty
{{\cos 2\pi k(x-a)} \over {k^{2n}}},$$
and
$$P_{2n-1}(x-a)=(-1)^n {{2(2n-1)!} \over {(2\pi)^{2n-1}}}\sum_{k=1}^\infty
{{\sin 2\pi k(x-a)} \over {k^{2n-1}}}.$$
We therefore obtain
$$P_{2n}(x-a)=(-1)^{n-1} {{2(2n)!} \over {(2\pi)^{2n}}}\sum_{k=1}^\infty
{1 \over {k^{2n}}}[\cos 2\pi k x+2\pi ka\sin 2\pi kx +O(a^2)],$$
and
$$P_{2n-1}(x-a)=(-1)^n {{2(2n-1)!} \over {(2\pi)^{2n-1}}}\sum_{k=1}^\infty
{1 \over {k^{2n-1}}}[\sin 2\pi kx-2\pi k a\cos 2\pi kx +O(a^2)].$$
These forms are then inserted into (2.8).  Noting that $C_k(1)=\gamma_k$, the $a^0$
term produces $\gamma_k$, $C_k(a) \to \gamma_1$ as $a \to 0$ and hence the result.
\qed

\medskip
\centerline{\bf Discussion}
\medskip

Here first we discuss the equivalence of Proposition 3 as a case of an addition formula which
we have previously presented \cite{coffeyseries} (Proposition 1).  We then show 
applications of differences of Stieltjes constants to some classic integrals of analytic
number theory.  We exhibit a new proof technique for certain log-log integrals.

As regards the Truesdell representation of $\zeta(s,kz)$, we note 
$${{s+n-1} \choose n}={{(-1)^n} \over {n!}}(1-s-n)_n={{(-1)^n} \over {n!}}{{\Gamma(1-s)}
\over {\Gamma(1-s-n)}}={1 \over {n!}}{{\Gamma(s+n)} \over {\Gamma(s)}}={1 \over {n!}}(s)_n,$$
so that
$$\zeta(s,kz)=\sum_{n=0}^\infty {{(s)_n} \over {n!}}(1-k)^n z^n  \zeta(s+n,z).$$
On the other hand, an old formula of Wilton \cite{wilton2} may be written as
$$\zeta(s,a+b)=\sum_{j=0}^\infty {{(-b)^j} \over {j!}}(s)_j\zeta(s+j,a), ~~|b|<|a|, ~~
\mbox{Re} ~a>0.$$
Thus the two formulas correspond with $b=-(1-k)z$ and $a=z$.  Lemma 1 of \cite{coffeyseries}
provides the derivative values
$$\left.\left({d \over {ds}}\right)^\ell (s)_j \right|_{s=1}=(-1)^{j+\ell} \ell! s(j+1,\ell+1),  $$
where $s(k,n)$ are the Stirling numbers of the first kind.  Therefore, from Proposition 1 of
\cite{coffeyseries}, we know the general form of the multiplication formula for the Stieltjes
constants,
$$\gamma_\ell(kz)=\gamma_\ell(z)+(-1)^\ell \sum_{j=2}^\infty {(k-1)^{j-1}z^{j-1} \over {(j-1)!}}
\sum_{k=0}^\ell (-1)^k {\ell \choose k} s(j,k+1) k!\zeta^{(\ell-k)}(j,z).$$
The Stirling numbers of the first kind may indeed be written with the generalized harmonic
numbers, and the first few are given by $s(n+1,1)=(-1)^n n!$, $s(n+1,2)=(-1)^{n+1}n!H_n$,
$s(n+1,3)=(-1)^n {{n!} \over 2}[H_n^2-H_n^{(2)}]$, and $s(n+1,4)=(-1)^{n+1}{{n!} \over 6}
[H_n^3-3H_nH_n^{(2)}+2H_n^{(3)}]$.    

We demonstrate how differences of Stieltjes constants may be used to efficiently evaluate some example log-log integrals, including 
$$I_{\pm}\equiv \int_0^1 {{\ln(-\ln x)} \over {1\pm x+x^2}}dx.$$
For $I_-$,
$$I_-=\int_0^1 \left({{1+x} \over {1+x^3}}\right)\ln(-\ln x)dx$$
$$=\sum_{m=0}^\infty (-1)^m\int_0^1 (1+x)x^{3m}\ln(-\ln x)dx$$
$$=\sum_{m=0}^\infty (-1)^m \int_0^\infty (1+e^{-u})e^{-(3m+1)u}\ln u ~du$$
$$=-\sum_{m=0}^\infty (-1)^m\left[\gamma\left({1 \over {3m+2}}+{1 \over {3m+1}}\right)
+{{\ln(3m+2)} \over {3m+2}}+{{\ln(3m+1)} \over {3m+1}}\right]$$
$$=-\gamma {{2\pi} \over {3\sqrt{3}}}-{1 \over 6}\sum_{m=0}^\infty \left[\ln 6\left({1 \over {m+1/3}}+{1 \over {m+1/6}}-{1 \over {m+5/6}}-{1 \over {m+2/3}}\right) \right.$$
$$\left. +{{\ln(m+1/3)} \over {m+1/3}}+{{\ln(m+1/6)} \over {m+1/6}}-{{\ln(m+5/6)} \over {m+5/6}}-{{\ln(m+2/3)} \over {m+2/3}}\right]$$
$$=-(\gamma+\ln 6){{2\pi} \over {3\sqrt{3}}}+{1 \over 6}\left[\gamma_1\left({2 \over 3}
\right)-\gamma_1\left({1 \over 3}\right)+\gamma_1\left({5 \over 6}\right)-\gamma_1\left(
{1 \over 6}\right)\right].$$
We used polygamma function values and (1.2).

For $I_+$,
$$I_+=\int_0^1 \left({{1-x} \over {1-x^3}}\right)\ln(-\ln x)dx$$
$$=\sum_{m=0}^\infty \int_0^1 (1-x)x^{3m}\ln(-\ln x)dx$$
$$=\sum_{m=0}^\infty \int_0^\infty (1-e^{-u})e^{-(3m+1)u}\ln u ~du$$
$$=\sum_{m=0}^\infty \left[\gamma\left({1 \over {3m+2}}-{1 \over {3m+1}}\right)
+{{\ln(3m+2)} \over {3m+2}}-{{\ln(3m+1)} \over {3m+1}}\right]$$
$$=-{{\gamma \pi} \over {3\sqrt{3}}}+{1 \over 3}\sum_{m=0}^\infty \left[{{\ln 3(m+2/3)}
\over {m+2/3}}-{{\ln 3(m+1/3)}\over {m+1/3}}\right]$$
$$=-(\gamma+\ln 3){\pi \over {3\sqrt{3}}}+{1 \over 3}\left[\gamma_1\left({2 \over 3}
\right)-\gamma_1\left({1 \over 3}\right)\right].$$

More generally,
$$I_{+n}=\int_0^1 {{\ln(-\ln x)} \over {x^{n-1}+x^{n-2}+\cdots+x+1}}dx$$
$$={{(\gamma+\ln n)} \over n}\left[\psi\left({1 \over n}\right)-\psi\left({2 \over n}\right)
\right]+{1 \over n}\left[\gamma_1\left({2 \over n}
\right)-\gamma_1\left({1 \over n}\right)\right],$$
$$I_{+n}^q=\int_0^1 {{x^q\ln(-\ln x)} \over {x^{n-1}+x^{n-2}+\cdots+x+1}}dx$$
$$={{(\gamma+\ln n)} \over n}\left[\psi\left({{q+1} \over n}\right)-\psi\left({{q+2} \over n}\right)\right]+{1 \over n}\left[\gamma_1\left({{q+2} \over n}
\right)-\gamma_1\left({{q+1} \over n}\right)\right], ~~\mbox{Re} ~q>-1,$$
and for $n$ odd,
$$I_{-n}^q =\int_0^1 {{x^q \ln(-\ln x)} \over {x^{n-1}-x^{n-2}+\cdots-x+1}}dx$$
$$={{(\gamma+\ln n)} \over {2n}}\left[\psi\left({{q+2} \over {2n}}\right)+\psi\left({{q+1} \over {2n}}\right)-\psi\left({{n+q+2} \over {2n}}\right)-\psi\left({{n+q+1} \over {2n}}\right)\right]$$
$$+{1 \over {2n}}\left[\gamma_1\left({{n+q+2} \over {2n}}\right)+\gamma_1\left({{n+q+1} \over {2n}}\right)-\gamma_1\left({{q+2} \over {2n}}\right)-\gamma_1\left({{q+1} \over {2n}}\right)
\right], ~~\mbox{Re} ~q>-1.$$

Similarly,
$$\int_0^1 {{\ln(-\ln x)} \over {1+x^2}}dx=-(\gamma+\ln 4){\pi \over 4}+{1 \over 4}
\left[\gamma_1\left({3 \over 4}\right)-\gamma_1\left({1 \over 4}\right)\right],$$
$$J_p\equiv \int_0^1 {{\ln(-\ln x)} \over {1+x^p}}dx={1 \over 2}\int_0^\infty {{e^{-(1-p/2)u}
\ln u} \over {\cosh(pu/2)}}du$$
$$={{(\gamma+\ln 2p)} \over {2p}}\left[\psi\left({1 \over {2p}}\right)-\psi\left({{p+1} \over {2p}}
\right)\right]+{1 \over {2p}}\left[\gamma_1\left({{p+1} \over {2p}}\right)-\gamma_1\left({1 \over {2p}}
\right)\right], ~~\mbox{Re} ~p>0,$$
and
$$J_p^2\equiv \int_0^1 {{\ln^2(-\ln x)} \over {1+x^p}}dx={1 \over {2p}}\left[\gamma^2+{\pi^2
\over 6}+2\gamma \ln 2p+\ln^2 2p\right]\left[\psi\left({{p+1} \over {2p}}\right)-\psi\left({1 \over {2p}}\right)\right]$$
$$+{1 \over p}(\gamma+\ln 2p)\left[\gamma_1\left({1 \over {2p}}\right)-\gamma_1\left({{p+1} \over {2p}}\right)\right]+{1 \over {2p}}\left[\gamma_2\left({1 \over {2p}}\right)-\gamma_2\left({{p+1} \over {2p}}\right)\right],$$
with the limits
$$\lim_{p \to \infty}J_p=-\gamma, ~~~~~~\lim_{p \to \infty}J_p^2=\gamma^2+\zeta(2).$$
As they should be, these limits are consistent with Corollary 3 and the more general
Proposition 5.

From the $J_p$ evaluation follows the integral identity
$$\int_0^1 \left\{{{(\gamma+\ln 2p)} \over {2p}}\left[\psi\left({1 \over {2p}}\right) -\psi\left({{p+1} \over {2p}}\right)\right]+{1 \over {2p}}\left[\gamma_1\left({{p+1} \over {2p}}\right)-\gamma_1\left({1 \over {2p}}\right)\right]\right\}dp$$
$$=\int_0^1 {{[\ln(2x)-\ln(x+1)]} \over {\ln x}}\ln(-\ln x)dx.$$
An analogous result applies for
$$\int_0^1 J_p^2 dp=\int_0^1 {{[\ln(2x)-\ln(x+1)]} \over {\ln x}}\ln^2(-\ln x)dx.$$

As an extension of $J_p$, for Re $p>0$ and Re $q>-1$, we have
$$J_p^q\equiv \int_0^1 {x^q \over {1+x^p}}\ln(-\ln x)dx$$
$$={{(\gamma+\ln 2p)} \over {2p}}\left[\psi\left({{q+1} \over {2p}}\right)-\psi\left({{p+q+1} \over {2p}}\right)\right]+{1 \over {2p}}\left[\gamma_1\left({{p+q+1} \over {2p}}\right)-\gamma_1\left({{q+1} \over {2p}}\right)\right],$$
with
$$\lim_{p \to \infty}J_p^q =-{{\gamma+\ln(q+1)} \over {q+1}}.$$

Propositions 1 and 4 apply to all of these integrals.  As a brief example, one finds
$$I_2\equiv\int_0^1 {{\ln(-\ln x)} \over {1+x^2}}dx={\pi \over 4}\left[\ln {{8\pi \Gamma^2(3/4)} \over {\Gamma^2(1/4)}}-\ln 4\right]={\pi \over 2}\ln {{\sqrt{2\pi} \Gamma(3/4)} \over {\Gamma(1/4)}}.$$

The value of $I_2$ has been known for a long time, and it may of course be written in many
equivalent forms.  However, the following method of evaluation may be new.
{\newline \bf Demonstration 1}.  
$$I_2\equiv\int_0^1 {{\ln(-\ln x)} \over {1+x^2}}dx=\int_{\pi/4}^{\pi/2} \ln(\ln(\tan x)) dx
={\pi \over 4}\left[\ln\left({\pi \over 8} \right)-2\ln{{\Gamma(3/4)} \over {\Gamma(5/4)}}\right].$$

The method below applies to a large variety of integrals, enabling another determination of
differences of Stieltjes constants at rational arguments.  A key feature of these integrals
is integrands with polynomial denominators with zeros at roots of unity.

{\it Proof}.  Write
$$I_2={1 \over {2i}}\int_0^1 \left({1 \over {x-i}}-{1 \over {x+i}}\right)\ln(-\ln x) ~dx,$$
and then apply
$$\int_0^1 {{\ln(-\ln x)} \over {x-a}}dx=-\gamma \ln\left({{a-1} \over a}\right)-\left. {{\partial \mbox{Li}_s} \over {\partial s}}\right|_{s=1}(a^{-1}), \eqno(3.1)$$
and (2.6) for the partial derivative to find
$$I_2={{\gamma \pi} \over 4}+{\pi \over 2}\ln\left({\pi \over 2}\right)+{1 \over {2i}}
\sum_{n=1}^\infty {{\zeta'(1-n)} \over {n!}}\left({\pi \over 2}\right)^n i^n[1-(-1)^n]$$
$$=\pi\left[{{\gamma} \over 4}+{1 \over 2}\ln\left({\pi \over 2}\right)\right]
+\sum_{m=0}^\infty {{(-1)^m \zeta'(-2m)} \over {(2m+1)!}}\left({\pi \over 2}\right)^{2m+1}.$$
Next separate the $m=0$ term of the sum and use the functional equation of the zeta function,
$\pi^{1-z}\zeta(z)=2^z\Gamma(1-z)\zeta(1-z)\sin {{\pi z} \over 2}$, along with $\zeta(-2m)=0$
for $m\geq 1$, to determine that for $m \geq 1$, 
$2(-1)^m \zeta'(-2m)=(2m)!\zeta(2m+1)/(2\pi)^{2m}$.  There results
$$I_2={\pi \over 4}\left[\gamma+\ln\left({\pi \over 8}\right)+\sum_{m=1}^\infty {{\zeta(2m+1)}
\over {16^m(2m+1)}}\right].$$
Using (e.g., \cite{grad}, p. 939)
$$\sum_{n=1}^\infty {x^{2n+1} \over {2n+1}}\zeta(2n+1)=-\gamma x+{1 \over 2}\left[\ln \Gamma
(1-x)-\ln \Gamma(x+1)\right], ~~~~|x|<1, \eqno(3.2)$$
completes the evaluation \cite{ftnote}.  \qed

As a further indication of the applicability of this method, we mention
{\newline \bf Lemma 1}. For $-\pi <\delta \leq \pi$,
$$I_\omega \equiv \int_0^1 {{\ln(-\ln x)} \over {x^2-2x\cos \delta+1}}dx
=-{\pi \over {2\sin \delta}}\left[{\delta \over \pi}\ln(2\pi)-\ln \delta + \ln {{\Gamma
\left(1+{\delta \over {2p}}\right)} \over {\Gamma \left(1-{\delta \over {2p}}\right)}}\right].$$
We only sketch the proof, as this is a known integral.

{\it Proof}.  We let $\omega=e^{i\delta}$ and use the factorization 
$(x-\omega)(x-\omega^*)=x^2-(\omega+\omega^*)x+|\omega|^2=x^2-2x\cos \delta +1$, giving
$$I_\omega={1 \over {\omega-\omega^*}}\int_0^1\left({1 \over {x-\omega}}-{1 \over {x-\omega^*}}
\right)\ln(-\ln x) dx.$$
We employ the integral (3.1), the partial derivative (2.6), and finally the summation (3.2).
Along the way we use elementary relations such as $1/\omega^*=\omega$ and
$$\ln\left[{{(\omega^*-1)\omega} \over {\omega^*(\omega-1)}}\right]=\ln\left({{1-e^{i\delta}} \over {1-e^{-i\delta}}}\right)=\ln[e^{i(\delta+\pi)}]=i(\delta+\pi).$$ 
\qed

Similarly, integrals of the form
$$\int_0^1 {{\ln^2(-\ln x)} \over {p(x)}}dx$$
with $p$ a polynomial having as zeros roots of unity, may be evaluated with the aid of 
$$\int_0^1 {{\ln^2(-\ln x)} \over {x-a}}dx=\left(\gamma^2+{\pi^2 \over 6}\right) \ln\left(
{{a-1} \over a}\right)+2\gamma \left. {{\partial \mbox{Li}_s} \over {\partial s}}\right|_{s=1} (a^{-1})-\left.{{\partial^2 \mbox{Li}_s} \over {\partial s^2}}\right|_{s=1} (a^{-1}),$$
and the partial derivative (2.7).  In summary, this method is comprised of the use of
partial fractions, logarithmic differentiation of a polylogarithm integral, application of
the partial derivatives of Li$_s$ at $s=1$, application of the functional equation of the
Riemann zeta function, and summation to $\ln \Gamma$ constants where pertinent.

Many other integrals follow from the results of this paper.  For instance, for $|a|=1$
but $a\neq 1$ or Re $a>1$ we have
$$\int_0^1 {{\ln(-\ln x)} \over {(x-a)^{m+1}}}dx={\gamma \over m}(-1)^m\left[{1 \over {(a-1)^m}}
-{1 \over a^m}\right]-{1 \over {m!}}\left({\partial \over {\partial a}}\right)^m \left. {{\partial \mbox{Li}_s} \over {\partial s}}\right|_{s=1} \left({1 \over a}\right).$$
This follows from (3.1) and (2.6).  In particular,
$$\int_0^1 {{\ln(-\ln x)} \over {(x-a)^2}}dx={1 \over a}\left\{\gamma\left[{1 \over {\ln a}}
-{1 \over {(a-1)}}\right]+{{\ln(\ln a)} \over {\ln a}}+\sum_{n=0}^\infty {{\zeta'(-n)} \over
{n!}} \ln^n(a^{-1}). \right\}$$
Such a result may be combined with the use of partial fractions to yield yet other integrals.


\pagebreak


\begin{thebibliography}{99}
\bibitem{nbs}M. Abramowitz and I. A. Stegun,
{Handbook of Mathematical Functions, Washington, National Bureau of Standards
(1964).}
\bibitem{apostol}T. M. Apostol,
{Introduction to Analytic Number Theory, Springer Verlag, New York (1976);
corrected fourth printing, 1995.}
\bibitem{berndt}B. C. Berndt,
{On the Hurwitz zeta function, Rocky Mtn. J. Math. {\bf 2}, 151-157 (1972).}
\bibitem{blago}I. V. Blagouchine, 
{A theorem for the closed–-form evaluation of the first generalized Stieltjes constant at rational arguments, arXiv 1401.3724v1 (2014).}
\bibitem{briggs}W. E. Briggs, 
{Some constants associated with the Riemann zeta-function, 
Mich. Math. J. {\bf 3}, 117-121 (1955).}
\bibitem{coffeyjmaa}M. W. Coffey,
{New results on the Stieltjes constants:  Asymptotic and exact evaluation, 
J. Math. Anal. Appl. {\bf 317}, 603-612 (2006); arXiv:math-ph/0506061.}
\bibitem{coffeystdiffs}M. W. Coffey,
{On representations and differences of Stieltjes coefficients, and other relations,
arXiv/math-ph/0809.3277v2 (2008); to appear in Rocky Mtn. J. Math.}
\bibitem{coffeystv2}M. W. Coffey,
{The Stieltjes constants, their relation to the $\eta_j$ coefficients,
and representation of the Hurwitz zeta function, Analysis {\bf 30}, 383 (2010),
arXiv/math-ph/:0706.0343v2 (2007).}
\bibitem{coffeyseries}M. W. Coffey
{Series representations for the Stieltjes constants, to appear in
Rocky Mtn. J. Math., arxiv/math-ph/0905.1111 (2009).}
\bibitem{coffeydec11}M. W. Coffey,
{Certain logarithmic integrals, including solution of Monthly problem 11629, zeta
values, and expressions for the Stieltjes constants, arXiv:1201.3393 (2012).}
\bibitem{coffeyaddison}M. W. Coffey,
{Integral representations of functions and Addison-type series for mathematical constants,
arXiv:1006.2551 (2010)}
\bibitem{coffeyhypergeo}M. W. Coffey,
{Hypergeometric summation representations of the Stieltjes constants, Analysis {\bf 33}, 
121-142 (2013), arXiv:1106.5148 (2011).}
\bibitem{edwards}H. M. Edwards,
{Riemann's Zeta Function, Academic Press, New York (1974).}
\bibitem{erdelyi}A. Erd\'{e}lyi, W. Magnus, F. Oberhettinger, and G. G. Tricomi,
{Higher transcendental functions, Vol. I, McGraw-Hill (1953).}
\bibitem{fine}N. J. Fine,
{Note on the Hurwitz zeta-function, Proc. Amer. Math. Soc. {\bf 2}, 361-364 (1951).}
\bibitem{grad}I. S. Gradshteyn and I. M. Ryzhik,
{Table of Integrals, Series, and Products, Academic Press, New York (1980).}
\bibitem{hansen}E. R. Hansen and M. L. Patrick,
{Some relations and values for the generalized Riemann zeta function,
Math. Comp. {\bf 16}, 265-274 (1962).}
\bibitem{hardy}G. H. Hardy, 
{Note on Dr. Vacca's series for $\gamma$, Quart. J. Pure Appl. Math. {\bf 43},
215-216 (1912).}
\bibitem{hurwitz}A. Hurwitz,
{Einige Eigenschaften der Dirichlet'schen Functionen $F(s)=\sum \left({D \over n}\right)
{1 \over n^s}$, die bei der Bestimmung der Classenanzahlen bin\"{a}rer quadratischer Formen auftreten, Z. Math. Phys. {\bf XXXVII}, 86-101 (1882).}
\bibitem{ivic}A. Ivi\'{c}, 
{The Riemann Zeta-Function, Wiley New York (1985).}
\bibitem{karatsuba}A. A. Karatsuba and S. M. Voronin,
{The Riemann Zeta-Function, Walter de Gruyter, New York (1992).}
\bibitem{kluyver}J. C. Kluyver, 
{On certain series of Mr. Hardy, Quart. J. Pure Appl. Math. {\bf 50}, 185-192
(1927).}
\bibitem{knessl1}C. Knessl and M. W. Coffey,
{An effective asymptotic formula for the Stieltjes constants, Math. Comp. {\bf 80}, 379-386 (2011).}
\bibitem{knessl2}C. Knessl and M. W. Coffey,
{An asymptotic form for the Stieltjes constants $\gamma_k(a)$ and for a sum $S_\gamma(n)$ appearing under the Li criterion, Math. Comp. {\bf 80}, 2197-2217 (2011).}
\bibitem{kreminski}R. Kreminski,
{Newton-Cotes integration for approximating Stieltjes (generalized Euler)
constants, Math. Comp. {\bf 72}, 1379-1397 (2003).}  
\bibitem{mitrovic}D. Mitrovi\'{c},
{The signs of some constants associated with the Riemann zeta function,
Mich. Math. J. {\bf 9}, 395-397 (1962).} 
\bibitem{riemann}B. Riemann,
{\"{U}ber die Anzahl der Primzahlen unter einer gegebenen Gr\"{o}sse, 
Monats. Preuss. Akad. Wiss., 671 (1859-1860).}
\bibitem{stieltjes}T. J. Stieltjes,
{Correspondance d'Hermite et de Stieltjes, Volumes 1 and 2, Gauthier-Villars,
Paris (1905).}
\bibitem{titch}E. C. Titchmarsh,
{The Theory of the Riemann Zeta-Function, 2nd ed., Oxford University
Press, Oxford (1986).}
\bibitem{truesdell}C. Truesdell,
{On the addition and multiplication theorem of special functions, Proc. N. A. S.
{\bf 36}, 752-755 (1950).}
\bibitem{wilton2}J. R. Wilton, 
{A note on the coefficients in the expansion of $\zeta(s,x)$ in powers of 
$s-1$, Quart. J. Pure Appl. Math. {\bf 50}, 329-332 (1927).}
\bibitem{yue}N.-Y. Zhang and K. S. Williams,
{Some results on the generalized Stieltjes constants, Analysis {\bf 14}, 
147-162 (1994).}  
\bibitem{ftnote}
{Details of the evaluation of the integrals $I_\pm$ by this method are separately 
available from the author.}
\end{thebibliography}
\end{document}